\documentclass[12pt]{article}
\usepackage{amsmath, amsthm, amssymb}

\begin{document}

\title{Asymptotic Flatness of the Weil-Petersson Metric on 
Teichm\"{u}ller Space}
\author{Zheng Huang\footnote{This work was partially supported 
by a Nettie S. Autrey fellowship and NSF grants 9971563 and 
0139887.}}

\date{June 11, 2003}
\newtheorem{theorem}{Theorem}
\newtheorem*{thm1}{Theorem 1}

\newtheorem{cor}[theorem]{Corollary}
\newtheorem{lem}[theorem]{Lemma}
\newtheorem{pro}[theorem]{Proposition}
\newtheorem{rem}[theorem]{Remark}

\newcommand{\WP}{Weil-Petersson}
\newcommand{\TS}{Teichm\"{u}ller space}
\newcommand{\Tt}{Teichm\"{u}ller theory}
\newcommand{\hq}{holomorphic quadratic}
\newcommand{\RS}{Riemann surface}
\newcommand{\hm}{harmonic map}
\newcommand{\Sc}{sectional curvature}
\newcommand{\cd}{complex dimension}
\newcommand{\Bd}{Beltrami differential}
\newcommand{\im}{identity map}
\newcommand{\Hd}{Hopf differential}
\newcommand{\Ts}{tangent space}
\newcommand{\de}{differential equation}

\maketitle
\tableofcontents 
\section{Introduction}

Let $\Sigma$ be a smooth, closed {\RS} of genus $g$, with $n$ punctures 
and $3g - 3 + n > 1$. Let $M_{-1}$ denote the space of metrics of constant 
curvature $-1$ on $\Sigma$. The group, $Diff_{0}$, of diffeomorphisms 
homotopic to the identity acts by pull back on $M_{-1}$, and we can 
define {\TS} ${\mathcal {T}}_{g,n}$, to be the quotient space 
${M_{-1}} / Diff_{0}(\Sigma)$. Note that the constant $3g - 3 + n $ is 
fundamental in {\Tt}: it is the {\cd} of the {\TS} ${\mathcal {T}}_{g,n}$.

The {\WP} metric on {\TS} has been heavily studied. It is defined by the 
$L^2$ inner product on the space of {\hq} differentials, the cotangent 
space of {\TS}. This metric has many curious properties; for example, it 
is a Riemannian metric with negative {\Sc} (\cite {Wp86} \cite {T87}), but 
it is incomplete (\cite {Chu} \cite {Wp75}). Also this metric is 
K\"{a}hler (\cite {Ahl}), and there is an upper bound 
$\frac {-1}{2 \pi (g-1)}$ for the holomorphic {\Sc} and Ricci curvature 
(\cite {Wp86}), a result conjectured by Royden (\cite {Ro}). However, 
there is no known negative upper bound for the {\Sc}s.

In this paper, we investigate the asymptotic behavior of the {\Sc}s 
of the {\WP} metric on {\TS}. Our method is to investigate {\hm}s from a 
nearly noded surface to nearby hyperbolic structures. One of the 
difficulties in estimating {\WP} curvatures is working with the operator 
$D = -2(\Delta - 2)^{-1}$, which appears in Tromba-Wolpert's curvature 
formula. Our approach is to study the {\Hd}s associated to {\hm}s and the 
analytic formulas resulting from the harmonicity of the maps. There is a 
natural connection between the operator $D$ and the local variations of 
the energy of a {\hm} between surfaces (\cite {Wf89}). From this 
connection, we estimate the solutions to some ordinary {\de}s to derive 
our curvature estimates. Our estimates imply that even though the {\Sc}s 
are negative, they are not staying away from zero. More specifically, we 
prove our main theorem:

\begin{theorem}
If the {\cd} of {\TS} $\mathcal {T}$ is greater than 1, then the {\Sc} of 
the {\WP} metric is not pinched from above by any negative constant, i.e, there 
is no negative upper bound for the {\WP} {\Sc}.
\end{theorem}

We remark that this asymptotic flatness results from pinching two 
independent closed geodesics on the surface. We will discuss more asymptotic 
flatness in an upcoming paper. The proof of the main theorem provides a 
quantitative bound for the {\Sc}, namely, 
\begin{theorem}
If the {\cd} of {\TS} $\mathcal {T}$ is greater than 1, and let $l$ be the 
length of shortest geodesic along a path to a boundary point in {\TS}, then 
there exists a family of tangent planes with {\WP} {\Sc} of the order $O(l)$.
\end{theorem}

Recently Brock and Farb (\cite {BF}) showed that the {\WP} metric on 
$\mathcal{T}$ is Gromov hyperbolic if and only if 
$dim_{C}({\mathcal {T}})\leq 2$. They pointed out that if the {\WP} metric on 
$\mathcal{T}$ had curvature pinched from above by a negative constant, the 
Gromov hyperbolicity in the case of the surface being doubly-punctured torus 
or 5-punctured sphere would be an immediate consequence of the comparison 
theorems. In the end of their paper, they asked the following question: if 
$int(\Sigma)$ is homeomorphic to a doubly-punctured torus or 5-punctured 
sphere, are the {\Sc}s of the {\WP} metric bounded away from zero 
(\cite {BF})? Our paper gives a negative answer to this question.

Naturally associated to a {\hm} 
$w: (\Sigma,\sigma |dz|^2) \rightarrow (\Sigma,\rho |dw|^2)$ is a quadratic 
differential $\Phi (\sigma, \rho) dz^2$, which is holomorphic with respect to 
the conformal structure of $\sigma$. This association of a quadratic 
differential to a conformal structure then defines a map 
$\Phi: {\mathcal{T}} \rightarrow QD(\Sigma)$ 
from {\TS} $\mathcal{T}$ to the space of {\hq} differentials $QD(\Sigma)$. 
This map is known to be a homeomorphism (\cite {Wf89}).

As an important computational tool in geometry of {\Tt}, the method of {\hm}s 
has been studied by many people. In particular, the second variation of the 
energy of the {\hm} $w = w(\sigma,\rho)$ with respect to the domain structure 
$\sigma$ (or image structure $\rho$) at $\sigma = \rho$ yields the {\WP} 
metric on $\mathcal{T}$ (\cite {T87}, \cite {Wf89}), and one can also 
establish Tromba-Wolpert's curvature formula of the {\WP} metric from this 
method (\cite {J}, see also \cite {Wf89}).

The moduli space of {\RS}s admits a compactification, known as the 
Deligne-Mumford compactification (\cite {Mu}), and any element of the 
compactification divisor can be thought of as a {\RS} with nodes, a 
connected complex space where points have neighborhoods complex isomorphic to 
either $\{ |z| < \varepsilon \}$ (regular points) or 
$\{ zw = 0; |z|, |w| < \varepsilon \}$ (nodes). We can think of noded 
surfaces arising as elements of the compactification divisor through a 
pinching process: fix a family of simple closed curves on $\Sigma$ such that 
each component of the complement of the curves has negative Euler 
characteristic. The noded surface is topologically the result of identifying 
each curve to a point, the node (\cite {Be}).

While writing this paper, the author learned that Scott Wolpert 
(\cite {Wp02}) has obtained results related to ours.

The organization of this paper is as follows. In section 2 we give the 
necessary background, define our terms and introduce the notations. Section 
3 is devoted to proving our main theorem. The discussion of this purpose is 
broken into subsections: in $\S 3.1$, we collect the local variational 
formulas associated with the {\hm}s, hence connect the local variations of 
the energy of the {\hm} to the operator $D = -2(\Delta - 2)^{-1}$; in 
$\S 3.2$ we describe the so called ``model case" of the problem; namely, we 
pinch two core geodesics of two cylinders into two points, and study the 
asymptotic behavior of the {\hm}s between cylinders. In the model case of 
this problem, the surface is a pair of cylinders; in $\S 3.3$, we describe 
a family of asymptotically flat 2-dimensional subspaces. We also establish 
two ordinary {\de}s, and relate the solutions of the equations to the 
operator $D = -2(\Delta - 2)^{-1}$, which leads us to estimate the 
curvature. We will establish the estimates of terms in the curvature 
formula when the surface is a pair of cylinders; in $\S 3.4$ we construct 
families of maps which have small tension, and are close to the {\hm}s 
resulting from pinching process in $\S 3.3$; finally in $\S 3.5$ we prove 
our main theorem using the estimates in $\S 3.3$ and the construction in 
$\S 3.4$. In section 4, we comment on the case when the surface has 
finitely many punctures.

{\bf{Acknowlegement:}} This paper contains part of the author's doctoral 
dissertation at Rice University. It is his great pleasure to express his 
deepest gratitude to his thesis advisor, Michael Wolf, for suggesting this 
problem, many stimulating conversations and for his patience. He also wants 
to thank Professor Frank Jones for his generous help. 

\section{Notations and Background}  

Recall that $\Sigma$ is a fixed, oriented, smooth surface of genus 
$g \geq 1$, and $n \geq 0$ punctures where $3g - 3 + n > 1$. We denote 
hyperbolic metrics on $\Sigma$ by $\sigma |dz|^2$ and $\rho |dw|^2$, where 
$z$ and $w$ are conformal coordinates on $\Sigma $. By the uniformization 
theorem, the set of all similarly oriented hyperbolic structures $M_{-1}$ 
can be identified with the set of all conformal structures on $\Sigma$ with 
the given orientation. Equivalently, this is the same as the set of all 
complex structures on $\Sigma$ with the given orientation. And {\TS} 
$\mathcal{T}$ is defined to be 
\begin{center}
$\mathcal{T}$$ = {M_{-1}} / Diff_{0}(\Sigma)$
\end{center}

It is well-known that {\TS} $\mathcal{T}$ has a complex structure 
($\cite {Ahl}$) and the co{\Ts} at a point $\Sigma \in \mathcal{T}$ is the 
space of {\hq} differentials $\Phi dz^2$ on $\Sigma$.  On $\Sigma$, 
there is a natural pairing of quadratic differentials and {\Bd}s 
$\mu (z)\frac {d\bar{z}}{dz}$ given by 
\begin{center}
$<\mu,\Phi> = $ Re $\int_{\Sigma}\mu (z)\Phi (z)dzd\bar{z}$
\end{center}

The {\Ts} at $\Sigma$ is the space of {\Bd}s modulo the trivial ones which 
are ones such that $<\mu,\Phi> = 0$. The {\WP} metric on $\mathcal{T}$ is 
obtained by duality from the $L^2$-inner product on $QD(\Sigma)$
\begin{center}
$<\phi,\psi> = \int_{\Sigma} \frac {\phi \bar{\psi}}{\sigma}dzd\bar{z}$
\end{center}
where $\sigma dzd\bar{z}$ is the hyperbolic metric on $\Sigma$.

For a Lipschitz map 
$w: (\Sigma,\sigma |dz|^2) \rightarrow (\Sigma,\rho |dw|^2)$, we define 
the energy density of $w$ at a point to be 
\begin{center}
$e(w;\sigma,\rho)= \frac {\rho (w(z))}{\sigma (z)} |w_{z}|^2 + 
\frac {\rho (w(z))}{\sigma (z)} |w_{\bar{z}}|^2 $
\end{center}
and the total energy 
$ E(w;\sigma,\rho) = \int_{\Sigma} e(w;\sigma,\rho) \sigma dzd\bar{z}$.

A critical point of $E(w;\sigma,\rho)$ is called a {\hm}; it satisfies 
the Euler-Lagrange equation, namely, 
\begin{center} 
$w_{z \bar{z}} + \frac{\rho_w}{\rho} w_z w_{\bar{z}} = 0.$
\end{center}

The Euler-Lagrange equation for the energy is the condition for the 
vanishing of the {\em{tension}}, which is, in local coordinates,
\begin{center} 
$\tau (w) = \Delta w^{\gamma} + 
^{N} \Gamma_{\alpha \beta}^{\gamma}w_i^{\alpha}w_j^{\beta} = 0$
\end{center}

It is well-known (\cite {Al} \cite {ES} \cite {Ha} \cite {SY} 
\cite {Sa}) that given $\sigma, \rho$, there exists a unique {\hm} 
$w: (\Sigma,\sigma) \rightarrow (\Sigma,\rho)$ homotopic to the identity 
of $\Sigma$, and this map is in fact a diffeomorphism.

Thus we have a {\hq} differential (also called a {\Hd}) 
$\Phi dz^2 = \rho w_z {\bar{w}}_{z}dz^2$, and evidently
\begin{center}
$\Phi = 0 \Leftrightarrow w$ is conformal $\Leftrightarrow \sigma = \rho.$
\end{center}
where $\sigma = \rho$ means that $(\Sigma,\sigma)$ and $(\Sigma,\rho)$ are 
the same point in {\TS} $\mathcal{T}$. This describes $\Phi$ as a 
well-defined map $\Phi: \mathcal{T} \rightarrow 
$ $QD(\Sigma)$ from {\TS} $\mathcal{T}$ to the space of {\Hd}s $QD(\Sigma)$. 
This map is a homeomorphism (\cite {Wf89}). Also note that the map $w$ can 
be extended to surfaces with finitely many punctures.

We define two auxiliary functions as following:
\begin{center}
$\mathcal{H} = \mathcal{H}$$ (\sigma,\rho) = 
\frac {\rho (w(z))}{\sigma (z)} |w_{z}|^2$
\\
$\mathcal{L} = \mathcal{L}$$ (\sigma,\rho) = 
\frac {\rho (w(z))}{\sigma (z)} |w_{\bar{z}}|^2$
\end{center}

Hence the energy density is $e = \mathcal {H} + \mathcal {L}$. Many aspects 
connected with a {\hm} between {\RS}s can be written in terms of 
$\mathcal {H}$ and $\mathcal {L}$ (and $\Phi$). We denote on 
$(\Sigma,\sigma |dz|^2)$
\begin{center}
$\Delta = \frac {4}{\sigma} \frac {\partial^{2}}{\partial z \partial \bar{z}}, 
K(\rho) = 
- \frac {2}{\rho} \frac {\partial^{2}}{\partial w \partial \bar{w}} log\rho,
K(\sigma) = - \frac {2}{\sigma} \frac {\partial^{2}}
{\partial z \partial \bar{z}} log\sigma,$
\end{center}
where $\Delta$ is the Laplacian, and $K(\rho)$, $K(\sigma)$ are curvatures 
of the metrics $\rho$ and $\sigma$, respectively.

The Euler-Lagrange equation gives
\begin{center}
$\Delta log {\mathcal{H}} = -2K(\rho){\mathcal{H}} + 2K(\rho){\cal{L}} + 
2K(\sigma)$.
\end{center}

When we restrict ourselves to the situation $K(\sigma)=K(\rho)=-1$, we will have 
the following facts (\cite {Wf89}):
\begin{itemize}
\item
$\mathcal{H} >$ 0;
\item
The {\Bd} $\mu = \frac {w_{\bar{z}}}{w_z} = \frac {\bar{\Phi}}{\sigma \mathcal{H}}$;
\item
$\Delta log {\mathcal{H}} = 2{\mathcal{H}} - 2{\mathcal{L}} - 2$.
\end{itemize}

\section{{\WP} Sectional Curvature}

In this section, we will prove the main theorem of this paper, namely, 
\begin{thm1}
If the {\cd} of {\TS} $\mathcal {T}$ is 
greater than 1, then the {\Sc} of the {\WP} metric is not pinched from above 
by any negative constant, i.e, there is no negative upper bound for the 
{\WP} {\Sc}.
\end{thm1}

Combined with the fact that there is no lower bound for the {\Sc} 
(\cite {Sch}), We have

\begin{cor} 
There are no negative bounds for the {\Sc} 
of the {\WP} metric when $dim_{C}({\mathcal {T}}) > 1$.
\end{cor}

The rest of this paper will be devoted to proving our main theorem. The discussion 
will be broken into 5 subsections. In $\S 3.1$, we collect some variational 
formulas we will need in the rest of discussion; in $\S 3.2$, we mostly describe the 
``model case", i.e. we consider a family of pertubations of the {\im} 
between a pair of cylinders; in $\S 3.3$, we describe a family of subplanes of the 
{\Ts} of {\TS} ${\mathcal {T}}_{g,2}$, and estimate terms 
in the curvature formula when the surface is a pair of cylinders; in $\S 3.4$, we 
construct $C^{2,\alpha}$ maps between surfaces and show that the constructed maps have 
small tension and are close to the {\hm}s we get from the pinching process; in 
$\S 3.5$, we prove the main theorem.
\subsection {Local Variation}
We consider a family of {\hm}s $w(t)$ for $t$ small, where $w(0) =$ id, the 
{\im}. Denote by $\Phi(t)$ the family of {\Hd}s determined by $w(t)$. We can rewrite 
the equation $\Delta log {\mathcal{H}} = 2{\mathcal{H}} - 2{\mathcal{L}} - 2$ as 
\begin{center}
$\Delta log {\mathcal{H}}(t) = 2{\mathcal{H}}(t) - 
\frac {2|\Phi(t)|^2}{\sigma^2 {\mathcal{H}}(t)} - 2$
\end{center}

For this equation, we see that the maximum principle will force all the odd order 
$t$-derivatives of ${\mathcal{H}}(t)$ to vanish, since the above equation only depends 
on the modulus of $\Phi(t)$ and not on its argument(\cite {Wf89}), and 
${\mathcal{H}}(t)$ is real-analytic in $t$ (\cite {Wf91}).

Wolf computed the $t$-derivative of various geometric quantities associated with the 
{\hm}s, and we collect these local variational formulas into our next lemma:
\begin{lem} (\cite {Wf89})
For the above notations, we have
\begin{itemize}
\item
${\mathcal{H}}(t) \geq 1$, and ${\mathcal{H}}(t) \equiv 1 \Leftrightarrow t = 0$;
\item
${\dot{\mathcal{H}}}(t) = \partial / \partial t^{\alpha}|_{0}$
${\mathcal{H}}(t) = 0$;
\item
${\dot{\mu}} = \partial / \partial t^{\alpha}|_{0} \mu(t) = 
\bar{\Phi}_{\alpha}/{\sigma}$;
\item
${\ddot{\mathcal{H}}}(t) = 
\frac {\partial ^2}{\partial t^{\alpha} \partial t^{\bar\beta}}|_{0}
{\mathcal{H}}(t) = D \frac {\Phi_{\alpha} \Phi_{\bar{\beta}}}{\sigma^2}$, 
where $D = -2(\Delta-2)^{-1}$. Evidently $D$ is a self-adjoint compact integral 
operator with a positive kernel, and it is the identity on constant functions.
\end{itemize}
\end{lem}

With this lemma, we obtain a partial {\de} about 
$\ddot{\mathcal{H}}(t)$:
\begin{center}
$(\Delta -2)(\ddot{\mathcal{H}}(t)) = -2 \frac {\Phi_{\alpha} \Phi_{\bar{\beta}}}{\sigma^2}$
\end{center}
\subsection{Harmonic Maps Between Cylinders}

For the sake of simplicity of exposition, we assume that our surface has no punctures. We 
will comment on punctured case in section 4.

Before we jump into the proof of the main theorem, in this subsection, we consider the 
asymptotics of {\hm}s when our surface is a pair of cylinders. 
We consider the asymptotics of a family of pertubations of the {\im} between two 
cylinders, and we call this is the ``model case". In the model case, the surface is a 
pair of cylinders $M_0$ 
and $M_1$. In particular, consider the boundary value problem of harmonically mapping 
the cylinder
\begin{center} 
$M = [l^{-1}sin^{-1}(l), \pi l^{-1}-l^{-1}sin^{-1}(l)] \times [0,1] $ 
\end{center}
with boundary identification  
\begin{center}
$[{\frac {sin^{-1}(l)}{l}}, 
{\frac {\pi}{l}}-{\frac {sin^{-1}(l)}{l}}] \times \{ 0 \} = 
[{\frac {sin^{-1}(l)}{l}}, {\frac {\pi}{l}}-{\frac {sin^{-1}(l)}{l}}] \times \{ 1 \}$ 
\end{center}
where the hyperbolic length element on $M$ is $lcsc(lx)|dz|$, to the cylinder 
\begin{center} 
$N = [L^{-1}sin^{-1}(L), \pi L^{-1}-L^{-1}sin^{-1}(L)] \times [0,1] $ 
\end{center}
with boundary identification  
\begin{center}
$[{\frac {sin^{-1}(L)}{L}}, 
{\frac {\pi}{L}}-{\frac {sin^{-1}(L)}{L}}] \times \{ 0 \} = 
[{\frac {sin^{-1}(L)}{L}}, {\frac {\pi}{L}}-{\frac {sin^{-1}(L)}{L}}] \times \{ 1 \}$ 
\end{center}
where the hyperbolic length element on $N$ is $Lcsc(Lu)|dw|$. Here $l$ and $L$ are the 
lengths of the simple closed core geodesics in the corresponding cylinders.\\
Say $w = u + iv$ is this {\hm} between cylinders $M$ and $N$, where 
\begin{center}
$u(l,L;x,y) = u(l,L;x), v(x,y) = y .$ 
\end{center}
The Euler-Lagrange equation becomes 
\begin{center}
$u'' = Lcot(Lu)(u'^2-1)$.
\end{center}
with boundary conditions $u({\frac {sin^{-1}(l)}{l}}) = {\frac {sin^{-1}(L)}{L}}$ and 
$u({\frac {\pi}{2l}}) = {\frac {\pi}{2L}}$. Note that both $M$ and $N$ admit an 
anti-isometric reflection about the curves $\{ \frac {\pi}{2l} \} \times [0,1]$ and 
$\{\frac {\pi}{2L}\} \times [0,1]$.

Since the quadratic differential 
$\Phi = \rho w_z {\bar{w}}_{z} = {\frac {1}{4}} L^2 csc^{2}(Lu)(u'^2 -1)$ is 
holomorphic in $M_0$, we have 
\begin{center}
$0 = {\frac {\partial}{\partial \bar{z}}}(\rho^2 w_z {\bar{w}}_{z}) = 
{\frac {\partial}{\partial {x}}}({\frac {1}{8}}L^2 csc^{2}(Lu)(u'^2 -1))$
\end{center}

Therefore $L^2 csc^{2}(Lu)(u'^2 -1) = c_{0}(l,L)$, where $c_{0}(l,L)$ is independent of 
$x$, and $c_{0}(l,l) = 0$ since $u(x)$ is the {\im} when $L=l$.

So we have 
\begin{center}
$u' = \sqrt {1 + c_{0}(l,L)L^{-2} sin^2 (Lu)}$
\end{center}
with boundary conditions $u({\frac {sin^{-1}(l)}{l}}) = {\frac {sin^{-1}(L)}{L}}$ and 
$u({\frac {\pi}{2l}}) = {\frac {\pi}{2L}}$. So the solution to the Euler-Lagrange equation 
can be derived from the equation
\begin{center}
$\int_{L^{-1}sin^{-1}L}^u {\frac {dv}{\sqrt {1 + c_{0}(l,L)L^2 sin^2 (Lv)}}} = 
x - l^{-1}sin^{-1}(l)$
\end{center}
with $c_{0}(l,L)$ chosen such that 
$\int_{L^{-1}sin^{-1}L}^{\frac{\pi}{2L}}{\frac {dv}{\sqrt {1 + c_{0}(l,L)L^2 sin^2 (Lv)}}} 
= {\frac{\pi}{2l}} - l^{-1}sin^{-1}l$.

It is not hard to show that when $l \rightarrow 0$, the solution 
$u(l,L;x)$ converges to a solution $u(L;x)$ to the ``noded" problem, i.e., 
$M = [1,+\infty)$, where we require $u(L;1) = L^{-1}sin^{-1}(L)$ and 
$\lim_{x \rightarrow +\infty} {u(L;x)} = \frac {\pi}{2L}$.

This ``noded" problem has the explicit solution as following (\cite {Wf91})
\begin{center} 
$u(L;x) = L^{-1}sin^{-1}\{ \frac {1-{\frac{(1-L)}{(1+L)}e^{2L(1-x)}}}{1+
{\frac{(1-L)}{(1+L)}e^{2L(1-x)}}} \} $
\end{center}
with the holomorphic energy
\begin{center} 
$H_{0}(L;x) = {\frac{L^{2}x^2}{4}}[\frac {1+{\sqrt{\frac{(1-L)}{(1+L)}}}e^{L(1-x)}}
{1-{\sqrt{\frac{(1-L)}{(1+L)}}}e^{L(1-x)}}]^2$
\end{center}
\subsection{Asymptotically Flat Subplanes}

Consider the surface $\Sigma$ which is developing two nodes, i.e., we are pinching two 
nonhomotopic closed geodesics $\gamma_0$ and $\gamma_1$ on $\Sigma$ to two points, say 
$p_0$ and $p_1$. We denote $M_0$ and $M_1$ their pinching neighborhoods, i.e., two 
cylinders described as $M$ above centered at $\gamma_0$ and $\gamma_1$, respectively.

We define $M(l_0, l_1)$ be the surface with two of the Fenchel-Nielsen coordinates, 
namely, the hyperbolic lengths of $\gamma_0$ and $\gamma_1$, are $l_0$ and $l_1$, 
respectively.  When we set the length of these two geodesics equal to $l$ simultaneously, 
we will have a point $M(l) = M(l,l)$ in {\TS} ${\mathcal{T}}_{g}$. 
Note that as $l$ tends to zero, the surface is developing two nodes. At this point $M(l)$, 
we look at two directions. First, we fix $\gamma_1$ in $M_1$ having length $l$, and pinch 
$\gamma_0$ in $M_0$ into length $L = L(t)$, where $L(0) = l$. So the $t$-derivative of 
$\mu_0(t)$, the {\Bd} of the resulting {\hm}, at $t = 0$ represents 
a tangent vector, $\dot{\mu}_0$, of {\TS} ${\mathcal{T}}_{g}$ at $M(l)$; 
we denote the resulting {\hm} by $W_{0}(t): M(l,l) \rightarrow M(L(t),l)$. Then we 
fix $\gamma_0$ in $M_0$ having length $l$, and pinch $\gamma_1$ in $M_1$ into length 
$L = L(t)$, so the $t$-derivative of $\mu_1(t)$ at $t = 0$ represents another tangent vector, 
$\dot{\mu}_1$, at $M(l)$; we denote the resulting {\hm} by 
$W_{1}(t): M(l,l) \rightarrow M(l;L(t))$. These two tangent vectors $\dot{\mu}_0$ and 
$\dot{\mu}_1$ will span a 2-dimensional subspace of the {\Ts} 
$T_{M(l)}{\mathcal{T}}_g$ to ${\mathcal{T}}_g$, hence we obtain a family, $\Omega_l$, of 
2-dimensional subspaces of the {\Ts} of ${\mathcal{T}}_{g}$.

\begin{pro}
The {\WP} {\Sc}s of $\Omega_l$ tend to zero as $l$ tends to zero.
\end{pro}

It is immediate that this proposition implies our main theorem.

To begin the proof of the proposition, we look at the curvature tensor of the {\WP} 
metric, which is given by (\cite {Wp86})
\begin{center}
$R_{\alpha {\bar{\beta}} \gamma {\bar{\delta}}} =  
(\int_{\Sigma}D(\dot{\mu_{\alpha}} \dot{\bar{\mu_{\beta}}})
\dot{\mu_{\gamma}}\dot{\bar{\mu_{\delta}}} dA) + 
(\int_{\Sigma}D(\dot{\mu_{\alpha}} \dot{\bar{\mu_{\delta}}})
\dot{\mu_{\gamma}}\dot{\bar{\mu_{\beta}}} dA)$
\end{center}
where $dA$ is the area element. Then the curvature of $\Omega$ is $R/{\Pi}$, where 
(\cite {Wp86})
\begin{center}
$R = R_{0 {\bar{1}} 0 {\bar{1}}} - R_{0 {\bar{1}} 1 {\bar{0}}} - 
R_{1 {\bar{0}} 0 {\bar{1}}} + R_{1 {\bar{0}} 1 {\bar{0}}}$
\end{center}
and 
\begin{eqnarray*}
{\Pi} & = & 4 <\dot{\mu}_0,\dot{\mu}_0><\dot{\mu}_1,\dot{\mu}_1> - 
2 |<\dot{\mu}_0,\dot{\mu}_1>|^2 - 2 Re(<\dot{\mu}_0,\dot{\mu}_1>)^2 \nonumber \\
& = & 4 <\dot{\mu}_0,\dot{\mu}_0><\dot{\mu}_1,\dot{\mu}_1> - 
4 |<\dot{\mu}_0,\dot{\mu}_1>|^2
\end{eqnarray*}

In the rest of this subsection, we will estimate terms in the above curvature formula in the 
model case, i.e., the surface is a pair of cylinders $M_0$ and $M_1$.

We denote $\phi_0(t)$ as the {\Hd} corresponding to the cylinder map $w_{0}(t)$ 
in $M_0$, and $\phi_1(t)$ as the {\Hd} corresponding to the cylinder map 
$w_{1}(t)$ in $M_1$. Here $w_{0}(t): M_0(l) \rightarrow M_0(L(t))$ and 
$w_{1}(t): M_1(l) \rightarrow M_1(L(t))$ are {\hm}s described as $w = u(x)+iy$ in the 
last subsection. In $M_1$, we still use $\phi_0$ to denote the {\Hd} corresponding 
to {\hm} $W_{0}(t)$, while in $M_0$, we still use $\phi_1$ as the {\Hd} 
corresponding to {\hm} $W_{1}(t)$. We also denote $\mu_0$ and $\mu_1$ as the 
corresponding {\Bd}s to $\phi_0$ and $\phi_1$.

We denote $a = a(l) = l^{-1}sin^{-1}(l)$, and 
$b = b(l) = \pi l^{-1}-l^{-1}sin^{-1}(l)$. And in this paper, $A \sim B$ means 
$A/C < B < CA$ for some constant $C > 0$.

So in $M_0$, as in $\S{3.2}$, we choose $L(t)$ so that 
${\frac{d}{dt}}|_{t=0}c_{0}(t) = {\frac{d}{dL}}|_{L=l} c_{0}(l,L) = 4$. Here we 
recall that $c_{0}(l,L(t)) = L^{2}csc^{2}(Lu)(u'^{2} - 1)$. We notice that 
$\dot{c}_{0}$ is never zero for all $l > 0$, otherwise we would have 
$\dot{\bar {w}}_z = 0$ as $\dot{\phi} = \sigma \dot{\bar {w}}_z$. Hence $w$ is a 
constant map by rotational invariance of the map. Thus we have 
$\dot{\phi_0} = {\frac{d}{dt}}|_{t=0}({\frac{1}{4}}c_{0}(t)) = 1$ in $M_0$.

Also, we see that $|\dot{\phi_1}|{\arrowvert}_{M_0} = \zeta (x,l)$ for 
$x \in [a,b]$, where $\zeta (x,l)$ satisfies that $\zeta (x,l) \leq C_{1}x^{-4}$ for 
$x \in [a,\pi /2l]$, and $\zeta (x,l) \leq C_{1}(\pi/l-x)^{-4}$ for 
$x \in [\pi /2l,b]$, and $\zeta (x,0)$ decays exponentially in $[1, +\infty]$. Here 
$C_{1} = C_{1}(l)$ is positive and bounded as $l \rightarrow 0$. To see this, notice that 
$\dot{\phi}_1$ is holomorphic and $|\dot{\phi}_1|$ is positive, so $log|\dot{\phi}_1|$ 
is harmonic in the cylinder $M_0$. Hence we can express $log|\dot{\phi}_1|$ in a Fourier 
series $\Sigma a_{n}(x)exp(-iny)$, and we compute 
$0 = \Delta log|\dot{\phi}_1| = \Sigma (a''_n - n^{2}a_n)exp(-iny)$. We will see, 
in $\S 3.4$, that $\dot{\phi}_1$ is close to 0 in $M_0$. Hence we conclude the 
properties $\zeta$ has. Similarly, we have $\dot{\phi_1}|_{M_1} = 1$ and 
$|\dot{\phi_0}||_{M_1} = \zeta (x,l)$ for $x \in [a,b]$.

Note that these estimates imply, informally, that most of the mass of $|\phi_0|$ resides 
in the thin part associated to $\gamma_0$, and most of the mass of $|\phi_1|$ resides in 
the thin part associated to $\gamma_1$.

Now in $M_0,$ the corresponding {\Bd} is
\begin{center}
$\dot{\mu}_0 = {\frac {d}{dt}}|_{t=0}({\frac {w_{\bar{z}}}{w_{z}}}) = 
\dot{\bar{\phi_0}}/{\sigma}$
\end{center}
and 
$|\dot{\mu}_0|^{2}|_{M_0} = 
|\dot{\bar{\phi_0}}/{\sigma}|^{2}|_{M_0} = l^{-4}sin^{4}(lx)$.

Also, in that same neighborhood $M_0$, 
\begin{center}
$\dot{\mu}_1 = \dot{\bar{\phi_1}}/{\sigma}$
\end{center}

Hence 
\begin{center}
$|\dot{\mu}_1|^{2}|_{M_0} = |\dot{\bar{\phi_1}}/{\sigma}|^{2}|_{M_0} 
= l^{-4}sin^{4}(lx)\zeta^2 (x,l) $. 
\end{center}

Similarly, in $M_1$, we have 
\begin{center}
$|\dot{\mu}_0|^{2}|_{M_1} = l^{-4}sin^{4}(lx)\zeta^2 (x,l)$ \\
$|\dot{\mu}_1|^{2}|_{M_1} = l^{-4}sin^{4}(lx) $
\end{center} 

\begin{lem}
$1/{\Pi} = O(l^3)$.
\end{lem}

We recall that ${\Pi} = 4 <\dot{\mu}_0,\dot{\mu}_0><\dot{\mu}_1,\dot{\mu}_1> - 
4 |<\dot{\mu}_0,\dot{\mu}_1>|^2$, and compute the asymptotics in $l$ of each term. 
Using $|\dot{\mu}_0|^{2}|_{M_0} = l^{-4}sin^{4}(lx)$, and noting that 
$a = a(l) = l^{-1}sin^{-1}(l)$, and $b = b(l) = \pi l^{-1}-l^{-1}sin^{-1}(l)$, we 
have that
\begin{eqnarray*}
<\dot{\mu}_0,\dot{\mu}_0>|_{M_{0}} & = & \int_{M_{0}}|\dot{\mu_0}|^2 \sigma dxdy \\
 & = & \int_{0}^{1}\int_{a}^{b}|\dot{\mu}_0|^2 \sigma dxdy \\
 & = & \int_{0}^{1}\int_{a}^{b}l^{-2}sin^{2}lx dxdy \\
 & = & \frac {\pi}{2} l^{-3} - l^{-3}sin^{-1}(l) \\
 & \sim &  l^{-3}
\end{eqnarray*}
And using $|\dot{\mu}_1|^{2}|_{M_0} = l^{-4}sin^{4}(lx)\zeta^2 (x,l)$, and 
$\zeta (x,l) \leq C_{1}x^{-4}$ for $x \in [a,\pi /2l]$, we have
\begin{eqnarray*}
<\dot{\mu_1},\dot{\mu}_1>|_{M_{0}} & = & \int_{M_{0}}|\dot{\mu}_1|^2 \sigma dxdy \\ 
 & = & 2\int_{0}^{1}\int_{a}^{\pi/2l}l^{-4}sin^{4}(lx)\zeta^2 (x,l) \sigma dxdy \\
 & \leq & 2C^{2}_{1}\int_{0}^{1}\int_{a}^{\pi/2l}l^{-2}sin^{2}(lx)x^{-8} dxdy \\
 & = &  O(1) 
\end{eqnarray*}

Also $<\dot{\mu}_0,\dot{\mu}_1> = \int_{\Sigma} \dot{\mu}_0 \dot{\bar{\mu}}_1 dA$, hence, 
\begin{eqnarray*}
<\dot{\mu}_0,\dot{\mu}_1>|_{M_{0}} & = & 
\int_{M_{0}}\dot{\mu}_0\dot{\mu}_1 \sigma dxdy\\
& \leq & C_{1}\int_{M_{0}}l^{-2}sin^{2}(lx)x^{-4} dxdy\\
& = & O(1)
\end{eqnarray*}

Similarly, 
\begin{center}
$<\dot{\mu}_0,\dot{\mu}_0>|_{M_{1}} = O(1) $\\
$<\dot{\mu}_1,\dot{\mu}_1>|_{M_{1}} \sim l^{-3} $\\
$<\dot{\mu}_0,\dot{\mu}_1>|_{M_{1}} = O(1)$
\end{center}

Note that ${\Pi} \geq (4<\dot{\mu}_0,\dot{\mu}_0><\dot{\mu}_1,\dot{\mu}_1> - 
4(<\dot{\mu}_0,\dot{\mu}_1>))^2|_{M_0} \sim l^{-3}$, which completes the proof of 
Lemma 6.

From Lemma 6, we have 
\begin{equation}
|R|/{\Pi} = O(|R|/(l^{-3})) = O(|R|l^3)
\end{equation}

Now we are left to estimate $|R|$. We have 

\begin{lem}
$|R| \leq 4 \int_{\Sigma}D(|\dot{\mu}_0|^2){|\dot{\mu}_1|}^{2}\sigma dxdy$ 
\end{lem}

Note that $D = -2(\Delta - 2)^{-1}$ is self-adjoint, hence we have 
\begin{center}
${\int_{\Sigma}}D(|\dot{\mu}_0|^2){|\dot{\mu}_1|}^{2}\sigma dxdy = 
{\int_{\Sigma}}D(|\dot{\mu}_1|^2){|\dot{\mu}_0|}^{2}\sigma dxdy$
\end{center}

Therefore,
\begin{eqnarray*}
R & = & R_{0 {\bar{1}} 0 {\bar{1}}} - R_{0 {\bar{1}} 1 {\bar{0}}} - 
R_{1 {\bar{0}} 0 {\bar{1}}} + R_{1 {\bar{0}} 1 {\bar{0}}}\\
& = & 
2 \int_{\Sigma}D(\dot{\mu}_0\dot{\bar{\mu}}_1)\dot{\mu}_0\dot{\bar{\mu}}_1 
\sigma dxdy + 
2 \int_{\Sigma}D(\dot{\mu_1}\dot{\bar{\mu}}_0)\dot{\mu_1}\dot{\bar{\mu}}_0 
\sigma dxdy\\
& - & \int_{\Sigma}D(|\dot{\mu_0}|^2){|\dot{\mu_1}|}^{2}\sigma dxdy - 
\int_{\Sigma}D(|\dot{\mu_1}|^2){|\dot{\mu_0}|}^{2}\sigma dxdy\\
& - & \int_{\Sigma}D(\dot{\mu_0}\dot{\bar{\mu}}_1)\dot{\mu_1}\dot{\bar{\mu}}_0 
\sigma dxdy - 
 \int_{\Sigma}D(\dot{\mu_1}\dot{\bar{\mu}}_0)\dot{\mu_0}\dot{\bar{\mu}}_1 
\sigma dxdy\\
& = & 2 \int_{\Sigma}D(\dot{\mu_0}\dot{\mu_1})\dot{\mu_1}\dot{\mu_0} \sigma dxdy - 
2 \int_{\Sigma}D(|\dot{\mu_0}|^2){|\dot{\mu_1}|}^{2}\sigma dxdy  
\end{eqnarray*}

The last equality follows from that here $\dot{\mu_0}$ and $\dot{\mu_1}$ are real 
functions. Now from lemma 4.3 of \cite {Wp86}, we have 
$|D(\dot{\mu_0}\dot{\mu_1})| \leq |D(|\dot{\mu_0}|^2)|^{1/2}|D(|\dot{\mu_1}|^2)|^{1/2}$. 
Then an application of H\"{o}lder inequality shows that 
\begin{eqnarray*}
|\int_{\Sigma}D(\dot{\mu_0}\dot{\mu_1})\dot{\mu_1}\dot{\mu_0} dA| & \leq & 
\int_{\Sigma}|D(|\dot{\mu_0}|^2)|^{1/2}|D(|\dot{\mu_1}|^2)|^{1/2}\dot{\mu_1}\dot{\mu_0}dA \\
& \leq & (\int_{\Sigma}D(|\dot{\mu_0}|^2){|\dot{\mu_1}|}^{2}dA)^{\frac {1}{2}}
(\int_{\Sigma}D(|\dot{\mu_1}|^2){|\dot{\mu_0}|}^{2}dA)^{\frac {1}{2}}\\
& = & \int_{\Sigma}D(|\dot{\mu_0}|^2){|\dot{\mu_1}|}^{2}\sigma dxdy 
\end{eqnarray*}

So we will have 
\begin{equation}
|R| \leq 4 \int_{\Sigma}D(|\dot{\mu_0}|^2){|\dot{\mu_1}|}^{2}\sigma dxdy
\end{equation}
which completes the proof of Lemma 7.

Therefore together with (1) and (2), to show the proposition, it is enough to 
show that 
$\int_{\Sigma}D(|\dot{\mu_0}|^2){|\dot{\mu_1}|}^{2}\sigma dxdy = o(l^{-3})$. 
In fact, 
\begin{lem}
$\int_{\Sigma}D(|\dot{\mu_0}|^2){|\dot{\mu_1}|}^{2}\sigma dxdy = O(l^{-2})$
\end{lem}

Firstly, from Lemma 4, we have
\begin{center}
$D(|\dot{\mu_0}|^2) = 
-2(\Delta-2)^{-1} \frac {|\dot{\phi_0}|^2}{\sigma^2} $
\end{center}

We recall in $M_0$, the {\Hd} $\phi_0$ is corresponding to the cylinder map 
$w_0(t): (M_0, \sigma) \rightarrow (M_0, \rho(t))$, and $|\dot{\phi_0}| = 1$. 
As in section 2 and $ \S 3.1$, we write 
${\mathcal{H}} = \frac {\rho (w(z))}{\sigma (z)} |w_{z}|^2$, therefore we can 
write $D(|\dot{\mu_0}|^2) = \ddot{\mathcal{H}}$. Then in the cylinder $M_0$, 
we have 
\begin{equation}
(\Delta-2) \ddot{\mathcal{H}} = -2 \frac {|\dot{\phi_0}|^2}{\sigma^2} = 
-2 l^{-4}sin^{4}(lx)
\end{equation}

A maximum principle argument implies that $\ddot{\mathcal{H}}$ is positive. 
This $\ddot{\mathcal{H}}$ converges to the holomorphic energy ${\ddot{H}_0}$ 
in the ``noded" problem (of $\S 3.2$) when $x$ is fixed but sufficiently 
large in $[a,\pi/2l]$. This convergence guarantees that $\ddot{\mathcal{H}}$ 
is bounded on the compacta in [a,b]. Thus we can assume that 
$A_{1}(l) = {\ddot{\mathcal{H}}}(a) = {\ddot{\mathcal{H}}}(l^{-1}sin^{-1}l) = O(1) > 0$. 
Then ${\ddot{\mathcal{H}}}(x)$ solves the following {\de}
\begin{equation}
(l^{-2}sin^{2}(lx)){\ddot{\mathcal{H}}}'' -2{\ddot{\mathcal{H}}} = -2l^{-4}sin^{4}(lx)   
\end{equation}
with the conditions
\begin{center}
${\ddot{\mathcal{H}}}(l^{-1}sin^{-1}l) = A_{1}(l), {\ddot{\mathcal{H}}}'(\pi/{2l}) = 0$
\end{center}

Recall from $\S 3.1$ that all the odd order $t$-derivatives of ${\mathcal{H}}(t)$ vanish.
Also notice that 
\begin{center}
$J(x) = \frac {sin^{2}(lx)}{2l^4}$
\end{center}
is a particular solution to equation (4). Hence we can check, by the method of reduction 
of the solutions, the general solution to equation (4) with the assigned conditons has 
the form 
\begin{center}
${\ddot{\mathcal{H}}}(l;x) = J(x) + A_{2}cot(lx) + A_{3}(1-lxcot(lx))$
\end{center}
where coefficients $A_{2} = A_{2}(l)$ and $A_{3} = A_{3}(l)$ are constants independent of 
$x$ and we can check, by substituting the solution into the assigned conditions, that they 
satisfy that
\begin{center}
$A_{2} = \frac {\pi}{2}A_{3} = O(l^{-1})$
\end{center}

Noticing that $\zeta(x,l) \leq C_{1}x^{-4}$ in $[a,\pi/2l]$, we compute the following
\begin{eqnarray*}
{\int_{M_0}}{\ddot{\mathcal{H}}}(x){|\dot{\mu_1}|}^{2}\sigma dxdy & = & 
{\int^{1}_{0}\int^{b}_{a}}{\ddot{\mathcal{H}}}(x)(l^{-2}sin^{2}(lx))\zeta^{2}(x,l)dxdy \\
& \leq & 2C^{2}_{1}({\int^{1}_{0}\int^{\pi/2l}_{a}}J(x)l^{-2}sin^{2}(lx)x^{-8}dxdy \\
& + & {\int^{1}_{0}\int^{\pi/2l}_{a}}A_{2}cot(lx)l^{-2}sin^{2}(lx)x^{-8}dxdy \\
& + & {\int^{1}_{0}\int^{\pi/2l}_{a}}A_{3}(1-lxcot(lx))l^{-2}sin^{2}(lx)x^{-8}dxdy) \\
& = & O(l^{-2}) + O(l^{-2}) + O(l^{-2}) \\
& = & O(l^{-2})
\end{eqnarray*}

Therefore 
\begin{eqnarray}
{\int_{M_0}}D(|\dot{\mu_0}|^2){|\dot{\mu_1}|}^{2}\sigma dxdy & = & 
{\int^{b}_{a}}{\ddot{\mathcal{H}}}{|\dot{\mu_1}|}^{2}\sigma dx \nonumber \\ 
& = & O(l^{-2}) 
\end{eqnarray}

Now let us look at the term $\int_{M_1}D(|\dot{\mu_0}|^{2})|\dot{\mu_1}|^{2}
\sigma dxdy$. 
Note that in $M_1$, which we identify with $[a,b] \times [0,1]$, we have 
$D(|\dot{\mu_0}|^{2}) = \ddot{\mathcal{H}}$, here $\mu_0(t)$ and 
$\mathcal{H}$ come from the {\hm} $W_0(t): M(l,l) \rightarrow M(L(t),l)$, where 
$|\dot{\mu_0}||_{M_1} = l^{-2}sin^{2}(lx)\zeta (x,l)$, and ${\mathcal{H}}(t)$ solves 
\begin{equation}
(l^{-2}sin^{2}(lx))\ddot{\mathcal{H}}'' -2\ddot{\mathcal{H}} = 
-2l^{-4}sin^{4}(lx)\zeta^{2}(x,l)
\end{equation}
with the conditions
\begin{center}
$\ddot{\mathcal{H}}(l^{-1}sin^{-1}l) = B_1(l)$, 
$\ddot{\mathcal{H}}(\pi/{l}-l^{-1}sin^{-1}l) = B_2(l)$.
\end{center}
Here $B_1(l)$ and $B_2(l)$ are positive and bounded as $l \rightarrow 0$, since $\ddot{\mathcal{H}}$ 
converges to the holomorphic energy in the ``noded" problem (of $\S 3.2$, when 
$M_1 = [1,\infty)$). We recall that $|\dot{\phi}_1||_{M_0} = \zeta (x,l)$, where 
$\zeta (x,l) \leq C_{1}x^{-4}$ for $x \in [a,\pi/2l]$ 
and $\zeta (x,l) \leq C_{1}(\pi/l-x)^{-4}$ for $x \in [\pi/2l,b]$.

Consider the equation 
\begin{equation}
(l^{-2}sin^{2}(lx))Y'' -2Y = 0 
\end{equation}
with the boundary conditions that satisfy 
\begin{center}
$Y(l^{-1}sin^{-1}l) = O(1)$ and $Y(\pi/{l}-l^{-1}sin^{-1}l) = O(1)$, as $l \rightarrow 0$.
\end{center}
We claim that there exists some $h = O(l)$ such that $\ddot{\mathcal{H}} - h$ is a 
supersolution to (7) for $x \in [l^{-1/4}, b-l^{-1/4}]$. To see this, we notice that 
$2|\dot {\phi}_0|^{2}|_{M_1}l^{-4}sin^{4}(lx)$ decays rapidly as 
$x \rightarrow \pi /2l$ for small $l$. So there is a positive number $h = O(l)$  
such that $2|\dot {\phi}_0|^{2}|_{M_1}l^{-4}sin^{4}(lx) < 2h $ for $x \in [l^{-1/4}, b-l^{-1/4}]$. 
Now for $x \in [l^{-1/4},b]$ we have
\begin{equation}
(l^{-2}sin^{2}(lx))(\ddot{\mathcal{H}}-h)'' -2(\ddot{\mathcal{H}}-h) = 
2h - 2|\dot {\phi}_0|^{2}|_{M_1}l^{-4}sin^{4}(lx) > 0
\end{equation}

Notice that for any constant $\lambda$, we have that if $Y(x)$ solves the equation (7), then 
so does $\lambda Y(x)$. So up to multiplying by a bounded constant, we have 
$Y|_{\partial M_1} > (\ddot{\mathcal{H}} - h)|_{\partial M_1}$. 
Hence $\ddot{\mathcal{H}} - h$ is a subsolution to (7) for $x \in [l^{-1/4}, b-l^{-1/4}]$. We 
can check the solutions to (7) have the form of 
\begin{center}
$Y(l;x) = B_{3}cot(lx) + B_{4}(1-lxcot(lx))$
\end{center} 
where constants $B_{3} = B_{3}(l)$ and $B_{4} = B_{4}(l)$ satisfy, from the boundary conditions 
for the equation (7), that
\begin{center}
$B_{3} = O(l)$ and $B_{4} = O(l)$
\end{center}

Therefore in $[l^{-1/4}, b-l^{-1/4}]$, We have $\ddot{\mathcal{H}} \leq h + Y(x)$. Now,
\begin{eqnarray*}
{\int_{M_1}}Y(x){|\dot{\mu_1}|}^{2}\sigma dxdy & = & 
{\int^{1}_{0}\int^{b}_{a}}Y(x)(l^{-2}sin^{2}(lx))dxdy \\
& = & {\int^{1}_{0}\int^{b}_{a}}B_{3}cot(lx)(l^{-2}sin^{2}(lx))dxdy \\ 
& + & {\int^{1}_{0}\int^{b}_{a}}B_{4}(1-lxcot(lx))(l^{-2}sin^{2}(lx))dxdy \\
& = & O(l^{-2}) + O(l^{-2}) \\
& = & O(l^{-2})
\end{eqnarray*}

Also for $x \in [a,l^{-1/4}]\cup [b-l^{-1/4},b]$, we apply maximum principle to equation (6) 
and find that 
\begin{center}
$\ddot{\mathcal{H}} \leq max(l^{-4}sin^{4}(lx_0)\zeta^2 (x_0,l), 
sup({\ddot{\mathcal{H}}}|_{{\partial}(([a,l^{-1/4}]{\cup}[b-l^{-1/4},b]){\times}[0,1])}))$
\end{center}
for some $x_0 \in [a,l^{-1/4}]$. Hence using the properties of $\zeta (x,l)$ in $[a,b]$, a direct 
computation shows that 
\begin{center}
${\int^{1}_{0}\int^{l^{-1/4}}_{a}}\ddot{\mathcal{H}}(x){|\dot{\mu_1}|}^{2}dA = 
O(l^{-3/4}) = o(l^{-1})$\\
${\int^{1}_{0}\int^{b}_{b-l^{-1/4}}}
\ddot{\mathcal{H}}(x){|\dot{\mu_1}|}^{2}dA = O(l^{-3/4}) = o(l^{-1})$
\end{center}
With this and using $\ddot{\mathcal{H}} \leq h + Y(x)$ for $x\in [l^{-1/4},b-l^{-1/4}]$, we have
\begin{eqnarray}
{\int_{M_1}}D(|\dot{\mu_0}|^{2}){|\dot{\mu_1}|}^{2}dA & = & 
{\int^{1}_{0}\int^{b}_{a}}\ddot{\mathcal{H}}(x){|\dot{\mu_1}|}^{2}dA \nonumber\\
& = & {\int^{1}_{0}\int^{b-l^{-\frac {1}{4}}}_{l^{-1/4}}}\ddot{\mathcal{H}}(x){|\dot{\mu_1}|}^{2}dA 
 + {\int^{1}_{0}\int^{l^{-\frac {1}{4}}}_{a}}\ddot{\mathcal{H}}(x){|\dot{\mu_1}|}^{2}dA 
\nonumber\\
& + & {\int^{1}_{0}\int^{b}_{b-l^{-1/4}}}\ddot{\mathcal{H}}(x){|\dot{\mu_1}|}^{2}dA \nonumber\\
& \leq & {\int^{1}_{0}\int^{b}_{a}} (Y(x) + h) {|\dot{\mu_1}|}^{2}\sigma dxdy 
 + o(l^{-1}) \nonumber\\
& = & O(l^{-2}) + O(l^{-2}) + o(l^{-1}) = O(l^{-2}). 
\end{eqnarray}

Now combine the estimates of Lemma 7, (9), (10), we will have 
\begin{eqnarray}
|R| & \leq & 4 \int_{\Sigma}D(|\dot{\mu_0}|^2){|\dot{\mu_1}|}^{2}\sigma dxdy 
\nonumber\\
& = & 4 (\int_{M_0}D(|\dot{\mu_0}|^2){|\dot{\mu_1}|}^{2}dA + 
\int_{M_1}D(|\dot{\mu_0}|^2){|\dot{\mu_1}|}^{2}dA  \nonumber\\
& = & O(l^{-2})
\end{eqnarray}

Now with Lemma 6 we have
\begin{center}
$|R|/{\Pi} = O(l^{-2})l^3 = O(l)$
\end{center}
namely, $|R|/{\Pi}$ tends to zero as we pinch the geodesics.
\subsection {Constructed Maps}

In last subsection, we estimated the terms in the curvature formula in the model 
case, i.e., the surface is a pair of cylinders. Essentially, in each pinching 
neighborhood $M_i$, we used the cylinder map $w_i$ instead of the {\hm} $W_i$ 
during the computation, for $i = 0, 1$. Now we are in the general setting, i.e., the 
surface $\Sigma$ develops two nodes. In this subsection, we will construct a family of 
maps $G_i$ to approximate the {\hm} $W_i$, and the essential parts of this family 
are the {\im} of the surface restricted in the non-cylindral part and the cylinder 
map in each pinching neighborhood. We will also show that this constructed family $G_i$ 
is reasonably close to the {\hm}s $W_i$, for $i=0, 1$; hence we can use the 
estimates we obtained in the previous subsection to the general situation. 

We recall some of the notations from previous subsections. We still set 
$M_0$ and $M_1$ to be the pinching neighborhoods of the nodes $p_0$ and $p_1$, 
respectively. Also $W_0(t)$ is the {\hm} corresponding to fixing $\gamma_1$ in 
$M_1$ having length $l$, and pinching $\gamma_0$ in $M_0$ into length 
$L = L(t)$, where $L(0) = l$. Similarly $W_{1}(t)$ is the {\hm} corresponding to 
fixing $\gamma_0$ in $M_0$ into length $l$, and pinching $\gamma_1$ in $M_1$ having length 
$L = L(t)$. Let $w_0(t)$ and $w_1(t)$ be cylinder maps in the model case, we want to 
show that $W_0(t)$ is close to $w_0(t)$ in $M_0$ and is close to {\im} outside of $M_0$.

We denote subsets 
$\Sigma_0 = \{p \in \Sigma: dist(p, \partial M_0) > 1\}$, and 
$\Sigma_1 = \{p \in \Sigma: dist(p, \partial M_1) > 1\}$. Define the 
$1$-tube of $\partial M_0$ as 
$B(\partial M_0, 1) = \{p \in \Sigma: dist(p, \partial M_0) \leq 1\}$, and 
the $1$-tube of $\partial M_1$ as 
$B(\partial M_1, 1) = \{p \in \Sigma: dist(p, \partial M_1) \leq 1\}$. We 
can construct a $C^{2,\alpha}$ map $G_0: \Sigma \rightarrow \Sigma$ such 
that 
\begin{center}
$G_0(p) = \left \{
\begin{array}{cc}
w_0(t)(p), & p \in M_0 \cap \Sigma_0\\
p, & p \in (\Sigma_0 \backslash M_0) \\
g_{t}(p), & p \in B(\partial M_0, 1)  
\end{array} \right. $
\end{center}
where $g_{t}(p)$ in $B(\partial M_0, 1)$ is constructed so that it satisfies:
\begin{itemize}
\item
$g_{t}(p) = p$ for $p \in \partial (\Sigma_0 \backslash (M_0 \cup B(\partial M_0,1)))$, 
and $g_{t}(p) = w_{0}(t)(p)$ for $p \in \partial (M_0 \cap \Sigma_0)$;
\item
$g_{t}$ is the {\im} when $t = 0$;
\item
$g_{t}$ is smooth and the tension of $g_t$ is of the order $O(t)$.
\end{itemize}
We note that $G_0$ consists of three parts. It is the cylinder map of $M_0$ in 
$M_0 \cap \Sigma_0$, the {\im} in $\Sigma_0 \backslash M_0$, and a smooth 
map in the intersection region $B(\partial M_0, 1)$. Among three parts of the 
constructed map $G_{0}(t)$, two of them, the {\im} and the cylinder map, 
are harmonic, hence have zero tension; thus the tension of $G_0(t)$ is 
concentrated in $B(\partial M_0, 1)$. From \S 3.2, for the cylinder map 
$w_0 = u(x) + iy$, we have $u' = \sqrt{1 + c_0(t)L^{-2}sin^{2}Lu}$, where 
$c_0(0) = 0, \dot{c}_0(0) = 4$,. Hence for 
$x \in [l^{-1}sin^{-1}(l), l^{-1}sin^{-1}(l) + 1]$,
\begin{center}
$w_{0,z}(x,y) = {\frac {1}{2}}(u'(x) + 1) 
= {\frac {1}{2}}((1 + O(1)c_{0}(t))^{\frac {1}{2}}+1) = 1 + O(1)t + O(t^2)$\\
$|w_{0,z}(x,y) - 1| = O(t) \rightarrow 0, (t \rightarrow 0)$\\
$|w_{0,z{\bar{z}}}(x,y)| = |{\frac {1}{4}}u''(x)| = O(|Lcot(Lu)(u'^2 - 1)|) = O(t)$
\end{center}

Thus we can require that $|g_{t,z}| \leq C_{2}t$ and 
$|g_{t,z\bar{z}}| \leq C_{2}t$, and the constant $C_2=C_{2}(t,l)$ is bounded 
in both $t$ and $l$ since the coefficient of $t$ for $c_0(t)$ is bounded for 
small $t$ and small $l$. With the local formula of the tension in section 2, 
we have $\tau(G_0(t))$, the tension of $G_0(t)$ is of the order $O(t)$.

Similarly, we have a $C^{2,\alpha}$ map $G_1: \Sigma \rightarrow \Sigma$. 
Note that these constructed maps $G_0$ and $G_1$ are not necessarily harmonic.

Now we are about to compare the constructed family $G_0(t)$ and the family of 
the harmonic maps $W_0(t)$. To do this, we consider the function 
$Q_0 = cosh(dist(W_0,G_0)) - 1$. 

\begin{lem}
$dist(W_0, G_0) \leq C_{3}t$ in $B(\partial M_0, 1)$, 
where the constant $C_3 = C_{3}(t,l)$ is bounded for small $t$ and $l$.
\end{lem}

First, we want to show that $Q_0$ is a $C^2$ function. Notice that both 
the {\hm} $W_0(t)$ and the constructed map $G_0(t)$ are the {\im} when 
$t = 0$, and both families vary smoothly without changing homotopy type 
in $t$ for sufficiently small $|t|$ (\cite {EL}). For all $l > 0$, and 
for any $\varepsilon > 0$, there exists a $\delta$ such that for 
$|t| < \delta$, we have $|W_0(t) - W_0(0)| < {\frac {\varepsilon}{2}}$ 
and $|G_0(t) - G_0(0)| < {\frac {\varepsilon}{2}}$. Therefore the 
triangular inequality implies that $|W_0(t) - G_0(t)| < \varepsilon$. 
Since $l$ is positive, the Collar Theorem (\cite {Bu}) implies that the 
surface has positive injectivity radius $r$ bounded below, and we choose 
our $\varepsilon << r$, so $Q_0$ is well defined and smooth.

We follow an argument in \cite {HW}. For any unit 
$v \in T^{1}(B(\partial M_0, 1))$, the map $G_0$ satifies the inequality 
$|\|dG_0(v)\| - 1| = O(t)$, hence $|dG_0(v)|^2 > 1 - \varepsilon_0$ where 
$\varepsilon_0 = O(t)$, then we find that for 
any $x \in \Sigma$,
\begin{eqnarray}
\Delta Q_0 & \geq & min\{|dG_0(v)|^2:dG_0(v) \bot \gamma_x \} Q \nonumber \\
 & - & <\tau (G_0), \bigtriangledown d(\bullet,W_0)|_{G_0(x)}>sinh(dist(W_0,G_0))
\end{eqnarray}
where $\gamma_x$ is the geodesic joining $G_0(x)$ to $W_0(x)$ with initial tangent 
vector $-\bigtriangledown d(\bullet,W_0)|_{G_0(x)}$ and terminal tangent vector 
$\bigtriangledown d(G_0(x),\bullet)|_{W_0(x)}$. If $G_0(t)$ does not coincide with 
$W_0(t)$ on $B(\partial M_0, 1)$, we must have all maxima of $Q_0(t)$ on the 
interior of $B(\partial M_0, 1)$, at any such maximum, we apply the inequality 
$|dG_0(v)|^2 > 1 - \varepsilon_0$ to (12) to find
\begin{eqnarray*}
0 \geq \Delta Q_0 & \geq & (1 - \varepsilon_0)Q_0 - 
\tau(G_0)(sinh(dist(W_0, G_0)) 
\end{eqnarray*}

so that at a maximum of $Q_0$, we have 
\begin{center}
$Q_0 \leq {\frac {\tau(G_0)sinhdist(W_0,G_0)}{(1 - \varepsilon_0)}}$
\end{center}

We notice that $Q_0$ is of the order $dist^{2}(W_0,G_0)$ and 
$sinhdist(W_0,G_0)$ is of the order $dist(W_0,G_0)$, this implies that 
$dist(W_0,G_0)$ is of the order $O(t)$ in $B(\partial M_0, 1)$, which 
completes the proof of Lemma 8. 

\begin{rem}
Lemma 8 implies that $Q_0(t)$ is of the order $O(t^2)$ in 
$B(\partial M_0, 1)$. 
\end{rem}
Note that $B(\partial M_0, 1)$ contains the boundary of the cylinder 
$M_0 \cap \Sigma_0$, which we identify with $[a+1,b-1] \times [0,1]$, 
where, again, $a = a(l) = l^{-1}sin^{-1}(l)$, 
and $b = b(l) = \pi l^{-1}-l^{-1}sin^{-1}(l)$. While in the cylinder 
$M_0 \cap \Sigma_0$, we have the inequality 
\begin{eqnarray*}
\Delta Q_0 & \geq & (1 - \varepsilon_0)Q_0 - \tau(G_0)(sinh(dist(W_0, G_0))\\
& = & (1 - \varepsilon_0)Q_0 - \tau(G_0)(tanh(dist(W_0, G_0))(1+Q_0) \\
& = & (1 - \varepsilon_0 - \tau(G_0))Q_0 - \tau(G_0)(tanh(dist(W_0, G_0)) \\
& \geq & 1/2 Q_0 - C_{4}t^2
\end{eqnarray*}
where the constant $C_4$ is bounded for small $t$ and $l$. Therefore we 
find that $Q_0(z,t)$ decays rapidly in $z=(x,y)$ for $x$ close enough to 
$\pi/2l$. Hence we can assume that $dist(W_0,G_0)$ is at most of order 
$C't$ in $[a+1,b-1]$, here $C'=C'(x,l)$ is no greater that $C_{5}x^{-2}$ 
for $x \in [a+1,\pi/2l]$, and no greater than $C_{5}(\pi/l-x)^{-2}$ for 
$x \in [\pi/2l,b-1]$, where $C_5$ is bounded for small $t$ and $l$. Both 
maps $W_0$ and $G_0$ are harmonic in $M_0 \cap \Sigma_0$, so they are 
also $C^1$ close (\cite {EL}), i.e. we have 
$|W_{0, \bar{z}} - G_{0, \bar{z}}| \leq C_{5}x^{-2}t$ for small $t$ and 
$l$, when $x \in [a+1,\pi/2l]$. Thus we see that 
$|{\dot{W}}_{0,\bar{z}} - {\dot{G}}_{0,\bar{z}}| = C_{5}x^{-2}$, for 
$x \in [a+1,\pi/2l]$. Also, 
$|{\dot{W}}_{0,\bar{z}} - {\dot{G}}_{0,\bar{z}}| = C_{5}(\pi/l-x)^{-2}$, 
for $x \in [\pi/2l,b-1]$.

As before we denote $\phi_0$ and $\phi_1$ as {\Hd}s corresponding to 
{\hm}s $W_{0}(t)$ and $W_{1}(t)$, respectively. We also denote $\mu_0$ and 
$\mu_1$ as the corresponding {\Bd}s in $M_0 \cap \Sigma_0$ and 
$M_1 \cap \Sigma_1$, respectively. Write 
$\phi_{G_0} = \rho G_{0,z}(t)\bar{G}_{0,z}(t)$ and 
$\phi_{G_1} = \rho G_{1,z}(t)\bar{G}_{1,z}(t)$. Notice that in $M_0 \cap \Sigma_0$, 
map $G_0$ is the cylinder map hence harmonic, so $\phi_{G_0}$ is the {\Hd} 
corresponding to $G_0$. When $t = 0$ we have $W_{0} = G_{0} =$ identity and 
$\rho = \sigma$, hence we can differentiate 
$\phi_{G_0} = \rho G_{0,z}(t)\bar{G}_{0,z}(t)$ in $t$ at $t = 0$, and find that 
$|\dot{\phi}_{G_0} - \dot{\phi}_{W_0}| = 
\sigma |{\dot{W}}_{0,\bar{z}} - {\dot{G}}_{0,\bar{z}}| \leq C_{5}l^{2}x^{-2}csc^{2}(lx) = O(1)$ for 
$x \in [a+1,\pi/2l]$, and 
$|\dot{\phi}_{G_0} - \dot{\phi}_{W_0}| \leq C_{5}l^{2}(\pi/l-x)^{-2}csc^{2}(lx) = O(1)$ for 
$x \in [\pi/2l,b-1]$. Therefore we have proved
\begin{lem}
$|\dot{\phi}_{G_0} - \dot{\phi}_{W_0}| = O(1)$ for 
$x \in [a+1,b-1]$.
\end{lem}

Similarly we have $|\dot{\phi_0} - \dot{\phi}_{G_0}||_{M_1 \cap \Sigma_1} = O(1)$, and 
$|\dot{\phi_1} - \dot{\phi}_{G_1}||_{M_0 \cap \Sigma_0} = O(1)$, also 
$|\dot{\phi_1} - \dot{\phi}_{G_1}||_{M_1 \cap \Sigma_1} =  O(1)$.
\subsection{Proof of the Main theorem}

In last subsection, we constructed families of maps and found that they are 
reasonably close to the corresponding families of {\hm}s between surfaces. 
Now we are ready to adapt the estimates in the model case to the general 
setting, and prove the proposition in this situation, which will imply our 
main theorem.

Recall from $\S 3.3$ that the {\Sc} is represented by $R/{\Pi}$. 
With Lemma 6 in the model case, the inequality 
${\Pi} \geq {\Pi}_{M_0 \cap \Sigma_0} \sim l^{-3}$ 
still holds. Hence $|R|/{\Pi} = O(l^{3})|R|$, so it will be sufficient to 
show that $|R| = o(l^{-3})$. From Lemma 7, we have 
$|R| \leq 4 \int_{\Sigma}D(|\dot{\mu_0}|^2){|\dot{\mu_1}|}^{2}\sigma dxdy$, 
so we will show that Lemma 8 still holds, which implies desired curvature 
estimate immediately.

Now we are about to estimate 
$\int_{\Sigma}D(|\dot{\mu_0}|^2){|\dot{\mu_1}|}^{2}\sigma dxdy$, which breaks 
into 3 integrals as follows:
\begin{eqnarray}
\int_{\Sigma}D(|\dot{\mu_0}|^2){|\dot{\mu_1}|}^{2}\sigma dxdy & = & 
\int_{M_0 \cap \Sigma_0}D(|\dot{\mu_0}|^2){|\dot{\mu_1}|}^{2}\sigma dxdy \nonumber \\
& + & 
\int_{M_1 \cap \Sigma_1}D(|\dot{\mu_0}|^2){|\dot{\mu_1}|}^{2}\sigma dxdy \nonumber\\ 
& + & \int_{K}D(|\dot{\mu_0}|^2){|\dot{\mu_1}|}^{2}\sigma dxdy 
\end{eqnarray}
where $K$ is the compact set disjoint from $(M_0 \cap \Sigma_0) \cup (M_1 \cap \Sigma_1)$. 

For the third integral, from the previous discussion, because of the convergence of the harmonic 
maps to the {\hm}s of ``noded" problem, we have both $|{\dot{\mu}_0}|$ and 
$|{\dot{\mu}_1}|$ are bounded. The maximum principle implies that 
$D(|{\dot{\mu}_0}|^2) = {\ddot{\mathcal {H}}} \leq sup \{|{\dot{\mu}_0}|^2 \}$. Note that 
$K$ is compact, hence we have the third integral is the order of $O(1)$.

From Lemma 4, we have $(\Delta-2) \ddot{\mathcal{H}} = 
-2 \frac {|\dot{\phi}_0|^2}{\sigma^2}$, so we can rewirte (12) as 
\begin{eqnarray}
\int_{\Sigma}D(|\dot{\mu_0}|^2){|\dot{\mu_1}|}^{2}\sigma dxdy & = & 
\int_{M_0 \cap \Sigma_0}\ddot{\mathcal{H}}{|\dot{\mu_1}|}^{2}\sigma dxdy  \nonumber \\
& + & 
\int_{M_1 \cap \Sigma_1}\ddot{\mathcal{H}}{|\dot{\mu_1}|}^{2}\sigma dxdy + O(1)
\end{eqnarray}

Now we will look at the first and the second integrals in (13). Recall that 
$(\Delta-2) \ddot{\mathcal{H}}^G = -2 \frac {|\dot{\phi}_0^G|^2}{\sigma^2}$, where 
${\mathcal{H}}^G$ is the holomorphic energy corresponding to the model case when the 
{\hm} is the cylinder map, and $\phi_0^G$ is the quadratic differential corresponding 
to the constructed map $G_0$. We also denote $\mu_0^G$ to be the {\Bd} 
corresponding to $\phi_0^G$. We assign similar meanings for $\phi_1^G$ and $\mu_1^G$. Also 
recall that $|\dot{\phi}_0 - \dot{\phi}_0^G| = O(1)$ and 
$|\dot{\phi}_1 - \dot{\phi}_1^G| = O(1)$ in 
$(M_0 \cap \Sigma_0) \cup (M_1 \cap \Sigma_1)$, here $O(1)$ is bounded in $l$ for small $l$. 
So we can set some $\lambda = O(1)$ (bounded in $l$ for small $l$) such that 
$|\dot{\phi}_0|^2 < \lambda^2 |\dot{\phi}_0^G|^2$ and at the 
boundary of $(M_0 \cap \Sigma_0) \cup (M_1 \cap \Sigma_1)$ satisfies 
$\ddot{\mathcal{H}} < \lambda \ddot{\mathcal{H}}^G$. For example, we can take 
$\lambda = 1 + max_{\partial {K}}({\frac {\ddot{\mathcal{H}}}{\ddot{\mathcal{H}}^G}}, 
{\frac {|\dot{\phi}_0|}{|\dot{\phi}_{0}^G|}})$, this $\lambda = O(1)$ because at 
$\partial K = \partial ((M_0 \cap \Sigma_0) \cup (M_1 \cap \Sigma_1))$, both 
$\ddot{\mathcal{H}}$, $\ddot{\mathcal{H}}^G$, and 
$\frac {|\dot{\phi}_0|}{|\dot{\phi}_{0}^G|}$ are bounded. Therefore,
\begin{center}
$(\Delta-2) (\ddot{\mathcal{H}} - \lambda \ddot{\mathcal{H}}^G) = 
2 {\frac {\lambda^2|\dot{\phi}_0^G|^2 - |\dot{\phi}_0|^2}{\sigma^2}} > 0$
\end{center}

So $(\ddot{\mathcal{H}} - \lambda \ddot{\mathcal{H}}^G)$ is a subsolution to the 
equation $(\Delta-2)Y = 0$ whose solutions have the form of 
$Y(l;x) = B_{3}cot(lx) + B_{4}(1-lxcot(lx))$, where constants $B_{3}$ and $B_{4}$ 
satisfy that $B_{3} = O(l)$ and $B_{4} = O(l)$. Hence in 
$(M_0 \cap \Sigma_0) \cup (M_1 \cap \Sigma_1)$, we have 
$\ddot{\mathcal{H}} \leq \lambda \ddot{\mathcal{H}}^M + B_{3}cot(lx) + B_{4}(1-lxcot(lx))$. 
We apply this into (13) and find that,
\begin{eqnarray*}
\int_{\Sigma}D(|\dot{\mu_0}|^2){|\dot{\mu_1}|}^2\sigma dxdy & = & 
\int_{(M_0 \cap \Sigma_0) \cup (M_1 \cap \Sigma_1)}
\ddot{\mathcal{H}}{|\dot{\mu_1}|}^2\sigma dxdy  + O(1) \nonumber \\
& \leq & \int_{(M_0 \cap \Sigma_0)}(\lambda \ddot{\mathcal{H}}^G 
+ Y(l;x))({|\dot{\mu}_1|}^2)\sigma dxdy \nonumber \\
 & + & \int_{(M_1 \cap \Sigma_1)}(\lambda \ddot{\mathcal{H}}^G 
+ Y(l;x))({|\dot{\mu}_1|}^2)\sigma dxdy + O(1) \nonumber \\
& \leq & \int_{M_0 \cap \Sigma_0}(\lambda \ddot{\mathcal{H}}^G 
+ Y(l;x))({2|\dot{\mu}_1|}^2+2|\dot{\mu}_1-\dot{\mu}_1^G|^2)dA \nonumber \\
 & + & \int_{M_1 \cap \Sigma_1}(\lambda \ddot{\mathcal{H}}^G 
+ Y(l;x))({2|\dot{\mu}_1|}^2+2|\dot{\mu}_1-\dot{\mu}_1^G|^2)dA \nonumber \\
 & + & O(1)
\end{eqnarray*}

Recalling the computation in $\S 3.3$, and $|\dot{\mu}_1-\dot{\mu}_1^G| = 
|\dot{\phi}_1-\dot{\phi}_1^G|/\sigma$, where 
$|\dot{\phi}_1-\dot{\phi}_1^G| \leq C_{5}l^{2}x^{-2}csc^{2}(lx)$ for $x \in [a+1,\pi /2l]$, 
we have the following:
\begin{eqnarray*}
\int_{M_0}(\lambda \ddot{\mathcal{H}}^G)|\dot{\mu}_1^G|^2 \sigma dxdy & = & O(l^{-2})\\
\int_{M_1}(\lambda \ddot{\mathcal{H}}^G)|\dot{\mu}_1^G|^2 \sigma dxdy & = & O(l^{-2})\\
\int_{M_0}Y(l,x)|\dot{\mu}_1^G|^2 \sigma dxdy & = & O(1)\\
\int_{M_1}Y(l,x)|\dot{\mu}_1^G|^2 \sigma dxdy & = & O(l^{-2})\\
\int_{M_0 \cap \Sigma_0}(\ddot{\mathcal{H}}^G + Y(l;x))(|\dot{\mu}_1-\dot{\mu}_1^G|^2)dA & = & o(l^{-2})\\
\int_{M_1 \cap \Sigma_1}(\ddot{\mathcal{H}}^G + Y(l;x))(|\dot{\mu}_1-\dot{\mu}_1^G|^2)dA & = & o(l^{-2})\\
\end{eqnarray*}

Therefore $\int_{\Sigma}D(|\dot{\mu_0}|^2){|\dot{\mu_1}|}^{2}\sigma dxdy = O(l^{-2})$, with 
Lemma 7, we have 
$|R| \leq 4 \int_{\Sigma}D(|\dot{\mu_0}|^2){|\dot{\mu_1}|}^{2}\sigma dxdy = O(l^{-2})$. 
Apply this into the curvature formula, we have $|R|/{\Pi} = O(l) \rightarrow 0 $ as 
$l \rightarrow 0 $, which completes the proof of proposition.

\section {Punctured Surface Case}
In previous sections, the argument on a compact surface is enough to prove our main 
theorem. However, we remark here that we allow the surface has finitely many punctures.

For the case of punctured surfaces, the existence of a harmonic diffeomorphism between 
punctured surfaces has been investigated by Wolf (\cite {Wf91}) and Lohkamp (\cite {L}). 
In particular, Lohkamp (\cite {L}) showed that a homeomorphism between punctured surfaces 
is homotopic to a unique harmonic diffeomorphism with finite energy, and the holomorphic 
quadratic differential corresponding to the {\hm} in the homotopy class of the 
identity is a bijection between {\TS} of punctured surfaces and the 
space of {\hq} differentials.

In this case, the set $\Sigma \backslash ((M_0 \cap \Sigma_0) \cup (M_1 \cap \Sigma_1))$ is no 
longer compact. Let $K_0$ be a compact surface with finitely many punctures, and $\{K_m\}$ be a compact 
exhaustion of $K_0$. We now estimate $\int_{K_0}D(|\dot{\mu_0}|^2){|\dot{\mu_1}|}^{2} dA$. Let 
$H(t)$ be the holomorphic energy corresponding to the {\hm} 
$w(t): (K_0, \sigma) \rightarrow (K_0, \rho (t))$, then $H(t)$ is bounded above and below, and has 
nodal limit 1 near the punctures (\cite {Wf91}), hence both $|\dot{\mu_0}|^2$ and $|\dot{\mu_1}|^2$ 
have the order $o(1)$ near the punctures. To see this, we consider $K_0$ as the union of $K_m$ and 
disjoint union of finitely many punctured disks, each equipped hyperbolic metric 
$\frac {|dz|^2} {z^{2}log^{2}z}$. Then $|\dot{\mu_0}| = O (|z|log^{2}z) \rightarrow 0$ as $z$ tends 
to the puncture, since the quadratic differential has a pole of at most the first order. A similar 
result holds for $|\dot{\mu_1}|$. We notice that $\partial K_0$ is the boundary of the cylinders, 
where the {\hm}s converge to a solution to the ``noded" problem as $l \rightarrow 0$, hence 
$D(|\dot{\mu_0}|^2) = {\ddot {H}}(t)$ is bounded on $\partial K_0$. Therefore we apply Omori-Yau 
maximum principle (\cite {O} \cite {Y}) to $(\Delta -2){\ddot {H}} = -2|\dot{\mu_0}|^2$ on $K_0$ and 
obtain that 
$sup(D(|\dot{\mu_0}|^2)) \leq max(sup(|\dot{\mu_0}|^2), max(D(|\dot{\mu_0}|^2))|_{\partial K_0}) = O(1)$. 
Hence we have 
\begin{eqnarray*}
\int_{K_0}D(|\dot{\mu_0}|^2){|\dot{\mu_1}|}^{2} dA & \leq & 
\int_{K_0} sup(|\dot{\mu_0}|^2){|\dot{\mu_1}|}^{2} dA \nonumber \\
& \leq & O(1)O(1)Vol(K_0) \nonumber \\
& = & O(1)
\end{eqnarray*}
In other words, our proof carries over to the punctured case, which completes the proof of our Main theorem.

\noindent
Current Address: \\
Zheng Huang\\
Department of Mathematics \\
The University of Oklahoma\\
Norman, OK 73019\\
E-mail: huang@math.ou.edu
\end{document}